\documentclass[a4paper,leqno,12pt]{article}

\usepackage{amsmath}
\usepackage{amssymb}
\usepackage{amsthm}

\begin{document}

\bibliographystyle{amsalpha}
\newtheorem{Theorem}{Theorem}[section]
\newtheorem{Lemma}{Lemma}[section]
\newtheorem{Remark}{Remark}[section]
\newtheorem{Corollary}{Corollary}[section]
\newtheorem{Proposition}{Proposition}[section]
\newtheorem{Example}{Example}[section]
\newtheorem{Definition}{Definition}[section]
\newtheorem{Problem}{Problem}[section]
\newtheorem{Proof}{Proof}[section]

\title{Symmetry and holomorphy of the third-order ordinary differential equation satisfied by the third Painlev\'e Hamiltonian \\}

\author{By\\
Yusuke Sasano \\
(The University of Tokyo, Japan)}
\maketitle


\begin{abstract} We study symmetry and holomorphy of the third-order ordinary differential equation satisfied by the third Painlev\'e Hamiltonian.

\textit{Key Words and Phrases.} B{\"a}cklund transformation, Birational transformation, Holomorphy condition, Painlev\'e equations.

2000 {\it Mathematics Subject Classification Numbers.} 34M55; 34M45; 58F05; 32S65.
\end{abstract}

\section{Introduction}

In this paper, we study symmetry and holomorphy of a 2-parameter family of third-order ordinary differential equation:
\begin{align}\label{eq:1}
\begin{split}
\frac{d^3u}{dt^3}=&\frac{1}{2t^4\left(\frac{du}{dt} \right)}\left\{ \left(t^2\frac{d^2u}{dt^2}-(2\alpha_1-1)u-\alpha_0 \right) \left(t^2 \frac{d^2u}{dt^2}+(2\alpha_1-1)u+\alpha_0 \right) \right\}\\
&-\frac{4\frac{du}{dt}\left(u+t\frac{du}{dt} \right)}{t}+\frac{4t+(2\alpha_1-1)^2}{2t^2}\frac{du}{dt}-\frac{2}{t}\frac{d^2u}{dt^2},
\end{split}
\end{align}
where $\alpha_0,\alpha_1$ are complex constants.

The Hamiltonian of the third Painlev\'e system (see \cite{7,9,10,12})
\begin{equation}
u:=H_{III}(q,p,t;\alpha_0,\alpha_1)=\frac{q^2p(p-1)+q\{(1-2\alpha_1)p-\alpha_0\}+tp}{t}
\end{equation}
satisfies the equation \eqref{eq:1}.

In \cite{5,6,12}, it is well-known that the Hamiltonian of the third Painlev\'e system satisfies the second-order differential equation of binormal form. Since we study its phase space, we need to remake a normal form. In order to consider the phase space for the equation \eqref{eq:1}, by resolving its accessible singular locus we transform the system of rational type into the system of polynomial type:
\begin{equation}\label{eq:6}
  \left\{
  \begin{aligned}
   \frac{dx}{dt} =&y,\\
   \frac{dy}{dt} =&yz+\frac{2\alpha_1-1}{t^2}x+\frac{\alpha_0}{t^2},\\
   \frac{dz}{dt} =&-\frac{1}{2}z^2-\frac{4}{t}x-4y-\frac{2}{t}z+\frac{(2\alpha_1-1)(2\alpha_1-3)}{2t^2}+\frac{2}{t}.
   \end{aligned}
  \right. 
\end{equation}

In this paper, we present symmetry and holomorphy of the system \eqref{eq:6} (see Theorems \ref{th:5},\ref{th:6},\ref{th:2} and Proposition \ref{pro:3}). At first, we make its phase space. These patching data suggest its B{\"a}cklund transformations. For example, it is well-known that the patching data
\begin{equation}\label{corre:1}
r:(X,Y,Z)=\left(\frac{1}{x},-(xy+\alpha)x,z \right) \quad (\alpha \in {\Bbb C})
\end{equation}
suggests the symmetry condition:
\begin{equation}\label{corre:2}
s:(x,y,z;\alpha) \rightarrow \left(x+\frac{\alpha}{y},y,z;-\alpha \right).
\end{equation}
We note that
\begin{align*}
dX \wedge dY \wedge dZ=dx \wedge dy \wedge dz,\\
ds(x) \wedge ds(y) \wedge ds(z)=dx \wedge dy \wedge dz.
\end{align*}
In the case of the Painlev\'e systems, we can recover all birational symmetries from these patching data by using the correspondence \eqref{corre:1},\eqref{corre:2} when $z=0$ (see \cite{7,13}).

\begin{center}
\begin{tabular}{|c|c|c|c|c|} \hline 
Systsm & Holomorphy conditions & Transition functions & Symmetry    \\ \hline 
$P_{III}$ & $D_6^{(1)}$-surface & preserving 2-form & $(2A_1)^{(1)}$  \\ \hline 
$P_{II}$ & $E_7^{(1)}$-surface & preserving 2-form & $A_1^{(1)}$  \\ \hline 
\eqref{eq:6} & Figure 1 & preserving 3-form & $<s_0,s_1,\pi>$ \\ \hline 
\end{tabular}
\end{center}
Considering certain reductions of the soliton equations, some symmetries are closed in the Painlev\'e systems. It is still an open question whether there exist some symmetries of the soliton equations that are closed in the equations satisfied by the Painlev\'e Hamiltonians. In this paper, we consider this problem from the different viewpoints of certain reductions of the soliton equations.

The phase space (see Figure 1) for the system \eqref{eq:6} is obtained by gluing three copies of $ (x_i,y_i,z_i,t) \in {\mathbb C}^3 \times ({\Bbb C}-\{0\})$ via the birational transformations preserving 3-form{\rm:\rm}
\begin{equation*}
dx_i \wedge dy_i \wedge dz_i=dx \wedge dy \wedge dz \quad (i=0,1).
\end{equation*}
These patching data $r_i \ (i=0,1)$ are slightly different from the type \eqref{corre:1}. However, for $r_i \ (i=0,1)$, we find its B{\"a}cklund transformations $s_0,s_1$ by modifying the correspondence \eqref{corre:1},\eqref{corre:2}, that is, the patching data
\begin{equation}\label{corre:3}
\tilde{r}:(X,Y,Z)=\left(\frac{1}{x},-(xy+g(z))x,z \right) \quad (g(z) \in {\Bbb C}(t)[z])
\end{equation}
suggests the symmetry condition:
\begin{equation}\label{corre:4}
\tilde{s}:(x,y,z) \rightarrow \left(x+\frac{g(z)}{y},y,z \right).
\end{equation}
These B{\"a}cklund transformations $s_0,s_1$ and $\pi$ generate an infinite group (see Theorem \ref{th:2}). In particular, the composition $s_1s_0$ becomes its translation (see Proposition \ref{pro:2}). These properties are different from symmetry and holomorphy of the third Painlev\'e system (see \cite{7,10,11,Sa,Tsuda}).

We also show that the system \eqref{eq:6} with $(\alpha_0,\alpha_1)=\left(0,\frac{3}{2} \right)$ admits a special solution solved by classical transcendental functions (see Section 5):
\begin{equation}
F(a;t)=\sum_{k=0}^{\infty} \frac{1}{(a)_k}\frac{t^k}{k!}.
\end{equation}
Here the symbol $(a)_k$ denotes $(a)_k:=a(a+1) \cdots (a+k-1)$.

This paper is organized as follows. In Section 2, after we make the equation \eqref{eq:1}, by resolving its accessible singular locus we make the system of polynomial type. In Section 3, we review the notion of accessible singularity and local index, we will prove Proposition \ref{pro:1}. We also study the accessible singularity and local index of the system \eqref{eq:6}. In Section 4, we will study its symmetry and holomorphy. In Section 5, we will study its special solutions.

\section{Differential system of Polynomial type}

At first, we make the equation \eqref{eq:1}.
After differentiating once, we obtain
\begin{equation}
\frac{dH_{III}}{dt}=-\frac{q(p-qp+qp^2-\alpha_0-2\alpha_1 p)}{t^2}.
\end{equation}
We can express $p$ by using the variables $H_{III},\frac{dH_{III}}{dt}$:
\begin{equation}\label{eq:3}
p=t\frac{dH_{III}}{dt}+H_{III}.
\end{equation}
After differentiating again, we obtain
\begin{equation}
\frac{d^2H_{III}}{dt^2}=-\frac{2p-2qp+2qp^2-\alpha_0-2\alpha_1 p-H_{III}+t\frac{dH_{III}}{dt}}{t^2}.
\end{equation}
We can express $q$ by using the variables $H_{III},\frac{dH_{III}}{dt},\frac{d^2H_{III}}{dt^2}$:
\begin{equation}\label{eq:4}
q=\frac{2t^2 \frac{dH_{III}}{dt}}{\alpha_0+(2\alpha_1-1)H_{III}+t(2\alpha_1+1)\frac{dH_{III}}{dt}+t^2\frac{d^2H_{III}}{dt^2}}.
\end{equation}
After differentiating, we obtain the equation \eqref{eq:1}.

In order to consider the phase space for the equation \eqref{eq:1}, we make the rational transformation.

\begin{Proposition}\label{pro:1}
The rational transformations
\begin{equation}\label{eq:5}
  \left\{
  \begin{aligned}
   x =&u,\\
   y =&\frac{du}{dt},\\
   z =&\frac{\frac{d^2u}{dt^2}-\frac{(2\alpha_1-1)u+\alpha_0}{t^2}}{\frac{du}{dt}}
   \end{aligned}
  \right. 
\end{equation}
take the system of rational type to the system \eqref{eq:6} of polynomial type.
\end{Proposition}
After we give the notion of accessible singularities and local index, we will prove Proposition \ref{pro:1} in next section.

\section{Accessible singularities }

Let us review the notion of accessible singularity. Let $B$ be a connected open domain in $\Bbb C$ and $\pi : {\mathcal W} \longrightarrow B$ a smooth proper holomorphic map. We assume that ${\mathcal H} \subset {\mathcal W}$ is a normal crossing divisor which is flat over $B$. Let us consider a rational vector field $\tilde v$ on $\mathcal W$ satisfying the condition
\begin{equation*}
\tilde v \in H^0({\mathcal W},\Theta_{\mathcal W}(-\log{\mathcal H})({\mathcal H})).
\end{equation*}
Fixing $t_0 \in B$ and $P \in {\mathcal W}_{t_0}$, we can take a local coordinate system $(x_1,\ldots ,x_n)$ of ${\mathcal W}_{t_0}$ centered at $P$ such that ${\mathcal H}_{\rm smooth \rm}$ can be defined by the local equation $x_1=0$.
Since $\tilde v \in H^0({\mathcal W},\Theta_{\mathcal W}(-\log{\mathcal H})({\mathcal H}))$, we can write down the vector field $\tilde v$ near $P=(0,\ldots ,0,t_0)$ as follows:
\begin{equation*}
\tilde v= \frac{\partial}{\partial t}+g_1 
\frac{\partial}{\partial x_1}+\frac{g_2}{x_1} 
\frac{\partial}{\partial x_2}+\cdots+\frac{g_n}{x_1} 
\frac{\partial}{\partial x_n}.
\end{equation*}
This vector field defines the following system of differential equations
\begin{equation}\label{39}
\frac{dx_1}{dt}=g_1(x_1,\ldots,x_n,t),\ \frac{dx_2}{dt}=\frac{g_2(x_1,\ldots,x_n,t)}{x_1},\cdots, \frac{dx_n}{dt}=\frac{g_n(x_1,\ldots,x_n,t)}{x_1}.
\end{equation}
Here $g_i(x_1,\ldots,x_n,t), \ i=1,2,\dots ,n,$ are holomorphic functions defined near $P=(0,\dots ,0,t_0).$

\begin{Definition}\label{Def}
With the above notation, assume that the rational vector field $\tilde v$ on $\mathcal W$ satisfies the condition
$$
(A) \quad \tilde v \in H^0({\mathcal W},\Theta_{\mathcal W}(-\log{\mathcal H})({\mathcal H})).
$$
We say that $\tilde v$ has an {\it accessible singularity} at $P=(0,\dots ,0,t_0)$ if
$$
x_1=0 \ {\rm and \rm} \ g_i(0,\ldots,0,t_0)=0 \ {\rm for \rm} \ {\rm every \rm} \ i, \ 2 \leq i \leq n.
$$
\end{Definition}

If $P \in {\mathcal H}_{{\rm smooth \rm}}$ is not an accessible singularity, all solutions of the ordinary differential equation passing through $P$ are vertical solutions, that is, the solutions are contained in the fiber ${\mathcal W}_{t_0}$ over $t=t_0$. If $P \in {\mathcal H}_{\rm smooth \rm}$ is an accessible singularity, there may be a solution of \eqref{39} which passes through $P$ and goes into the interior ${\mathcal W}-{\mathcal H}$ of ${\mathcal W}$.

Here we review the notion of {\it local index}. Let $v$ be an algebraic vector field with an accessible singular point $\overrightarrow{p}=(0,\ldots,0)$ and $(x_1,\ldots,x_n)$ be a coordinate system in a neighborhood centered at $\overrightarrow{p}$. Assume that the system associated with $v$ near $\overrightarrow{p}$ can be written as
\begin{align}\label{b}
\begin{split}
\frac{d}{dt}Q\begin{pmatrix}
             x_1 \\
             x_2 \\
             \vdots\\
             x_n
             \end{pmatrix}=&\frac{1}{x_1}\left\{Q\begin{bmatrix}
             a_1 & & & &  \\
             & a_2 & & & \\
             & & \ddots & & \\
             & & & & a_n
             \end{bmatrix}Q^{-1}{\cdot}Q\begin{pmatrix}
             x_1 \\
             x_2 \\
             \vdots\\
             x_n
             \end{pmatrix}+\begin{pmatrix}
             x_1h_1(x_1,\ldots,x_n,t) \\
             h_2(x_1,\ldots,x_n,t) \\
             \vdots\\
             h_n(x_1,\ldots,x_n,t)
             \end{pmatrix}\right\},\\
             &(h_i \in {\Bbb C}(t)[x_1,\ldots,x_n], \ Q \in GL(n,{\Bbb C}(t)), \ a_i \in {\Bbb C}(t))
             \end{split}
             \end{align}
where $h_1$ is a polynomial which vanishes at $\overrightarrow{p}$ and $h_i$, $i=2,3,\ldots,n$ are polynomials of order at least 2 in $x_1,x_2,\ldots,x_n$, We call ordered set of the eigenvalues $(a_1,a_2,\cdots,a_n)$ {\it local index} at $\overrightarrow{p}$.

We are interested in the case with local index
\begin{equation}\label{integer}
(1,a_2/a_1,\ldots,a_n/a_1) \in {\Bbb Z}^{n}.
\end{equation}
These properties suggest the possibilities that $a_1$ is the residue of the formal Laurent series:
\begin{equation}
y_1(t)=\frac{a_1}{(t-t_0)}+b_1+b_2(t-t_0)+\cdots+b_n(t-t_0)^{n-1}+\cdots \quad (b_i \in {\Bbb C}),
\end{equation}
and the ratio $(a_2/a_1,\ldots,a_n/a_1)$ is resonance data of the formal Laurent series of each $y_i(t) \ (i=2,\ldots,n)$, where $(y_1,\ldots,y_n)$ is original coordinate system satisfying $(x_1,\ldots,x_n)=(f_1(y_1,\ldots,y_n),\ldots,f_n(y_1,\ldots,y_n)) \ f_i(y_1,\ldots,y_n) \in {\Bbb C}(t)(y_1,\ldots,y_n)$.

If each component of $(1,a_2/a_1,\ldots,a_n/a_1)$ has the same sign, we may resolve the accessible singularity by blowing-up finitely many times. However, when different signs appear, we may need to both blow up and blow down.

The $\alpha$-test,
\begin{equation}
t=t_0+\alpha T, \quad x_i=\alpha X_i, \quad \alpha \rightarrow 0,
\end{equation}
yields the following reduced system:
\begin{align}
\begin{split}
\frac{d}{dT}P\begin{pmatrix}
             X_1 \\
             X_2 \\
             \vdots\\
             X_n
             \end{pmatrix}=\frac{1}{X_1}P\begin{bmatrix}
             a_1(t_0) & & & &  \\
             & a_2(t_0) & & & \\
             & & \ddots & & \\
             & & & & a_n(t_0)
             \end{bmatrix}P^{-1}{\cdot}P\begin{pmatrix}
             X_1 \\
             X_2 \\
             \vdots\\
             X_n
             \end{pmatrix}, \ P \in GL(n,{\Bbb C}).
             \end{split}
             \end{align}
From the conditions \eqref{integer}, it is easy to see that this system can be solved by rational functions.

{\bf Proof of Proposition \ref{pro:1}.} At first, setting
$$
p=u, \quad q=\frac{du}{dt}, \quad r=\frac{d^2u}{dt^2}.
$$
Next, we calculate the singular points when $q=0$. In this case, we solve the equation:
$$
\left(t^2r-(2\alpha_1-1)p-\alpha_0 \right) \left(t^2r+(2\alpha_1-1)p+\alpha_0 \right)=0.
$$
We obtain 2 solutions:
$$
r=\frac{(2\alpha_1-1)p+\alpha_0}{t^2} \quad \rm{ or \rm} \quad -\frac{(2\alpha_1-1)p+\alpha_0}{t^2}.
$$
Let us take the coordinate system $(P,Q,R)$ centered at the point $(P,Q,R)=(p,0,\frac{(2\alpha_1-1)p+\alpha_0}{t^2})$:
$$
P=p, \quad Q=q, \quad R=r-\frac{(2\alpha_1-1)p+\alpha_0}{t^2}.
$$
This system is rewritten as follows:
\begin{align*}
\frac{d}{dt}\begin{pmatrix}
             P \\
             Q \\
             R 
             \end{pmatrix}&=\frac{1}{Q}\left\{\begin{pmatrix}
             0 & 0 & 0 \\
             0 & \frac{\alpha_0}{t^2} & 0 \\
             0 & 0 & \frac{\alpha_0}{t^2}
             \end{pmatrix}\begin{pmatrix}
             P \\
             Q \\
             R 
             \end{pmatrix}+\cdots\right\}
             \end{align*}
satisfying \eqref{b}. In this case, the local index is $\left(0,\frac{\alpha_0}{t^2},\frac{\alpha_0}{t^2} \right)$. This suggests the possibilities that we will resolve this accessible singular locus by blowing-up once in the direction $R$.

Blowing up along the curve $\{(P,Q,R)=(P,0,0)\}${\rm : \rm}
$$
x=P, \quad y=Q, \quad z=\frac{R}{Q},
$$
then we obtain the system \eqref{eq:6}. \qed

In order to consider the phase spaces for the system \eqref{eq:6}, let us take the compactification $[z_0:z_1:z_2:z_3] \in {\mathbb P}^3$ of $(x,y,z) \in {\mathbb C}^3$ with the natural embedding
$$
(x,y,z)=(z_1/z_0,z_2/z_0,z_3/z_0).
$$
Moreover, we denote the boundary divisor in ${\mathbb P}^3$ by $ {\cal H}$. Extend the regular vector field on ${\mathbb C}^3$ to a rational vector field $\tilde v$ on ${\mathbb P}^3$. It is easy to see that ${\mathbb P}^3$ is covered by four copies of ${\mathbb C}^3${\rm : \rm}
\begin{align*}
&U_0={\mathbb C}^3 \ni (x,y,z),\\
&U_j={\mathbb C}^3 \ni (X_j,Y_j,Z_j) \ (j=1,2,3),
\end{align*}
via the following rational transformations
\begin{align*}
& X_1=1/x, \quad Y_1=y/x, \quad Z_1=z/x,\\
& X_2=x/y, \quad Y_2=1/y, \quad Z_2=z/y,\\
& X_3=x/z, \quad Y_3=y/z, \quad Z_3=1/z.
\end{align*}

\begin{Lemma}
The rational vector field $\tilde v$ has two accessible singular loci{\rm : \rm}
\begin{equation}
  \left\{
  \begin{aligned}
   C_1 &=\{(X_1,Y_1,Z_1)|X_1=Z_1=0\} \cup \{(X_2,Y_2,Z_2)|Y_2=Z_2=0\} \cong {\Bbb P}^1,\\
   P_2 &=\{(X_3,Y_3,Z_3)|X_3=Y_3=Z_3=0\},\\
   \end{aligned}
  \right. 
\end{equation}
where $C_1 \cong {\Bbb P}^1$ is multiple locus of order $2$.
\end{Lemma}

Next let us calculate its local index at the point $P_2$.
\begin{center}
\begin{tabular}{|c|c|c|} \hline 
Singular point & Type of local index   \\ \hline 
$P_2$ & $(\frac{1}{2},\frac{3}{2},\frac{1}{2})$  \\ \hline 
\end{tabular}
\end{center}

In order to do analysis for the accessible singular locus $C_1$, we need to replace a suitable coordinate system because each point has multiplicity of order 2.

At first, let us do the Painlev\'e test. To find the leading order behaviour of a singularity at $t=t_1$ one sets
\begin{equation*}
  \left\{
  \begin{aligned}
   x & \propto \frac{a}{(t-t_1)^m},\\
   y & \propto \frac{b}{(t-t_1)^n},\\
   z & \propto \frac{c}{(t-t_1)^p},
   \end{aligned}
  \right. 
\end{equation*}
from which it is easily deduced that
\begin{equation*}
m=1, \quad n=2, \quad p=1.
\end{equation*}
The order of pole $(m,n,p)=(1,2,1)$ suggests a suitable coordinate system to do analysis for the accessible singularities, which is explicitly given by
\begin{equation*}
(X^{(1)},Y^{(1)},Z^{(1)})=\left(\frac{x}{z},\frac{y}{z^2},\frac{1}{z} \right).
\end{equation*}

In this coordinate, the singular point is given as follows:
\begin{equation*}
   P= \left\{(X^{(1)},Y^{(1)},Z^{(1)})=\left(-\frac{1}{2},-\frac{1}{4},0 \right)\right\}.
\end{equation*}

Let us take the coordinate system $(p,q,r)$ centered at the point $P$:
$$
p=X^{(1)}+\frac{1}{2}, \quad q=Y^{(1)}+\frac{1}{4}, \quad r=Z^{(1)}.
$$
The system \eqref{eq:6} is rewritten as follows:
\begin{align*}
\frac{d}{dt}\begin{pmatrix}
             p \\
             q \\
             r 
             \end{pmatrix}&=\frac{1}{r}\left\{\begin{pmatrix}
             -\frac{1}{2} & -1 & 0 \\
             0 & -2 & 0 \\
             0 & 0 &  -\frac{1}{2}
             \end{pmatrix}\begin{pmatrix}
             p \\
             q \\
             r 
             \end{pmatrix}+\cdots\right\}
             \end{align*}
satisfying \eqref{b}. To the above system, we make the linear transformation
\begin{equation*}
\begin{pmatrix}
             X \\
             Y \\
             Z 
             \end{pmatrix}=\begin{pmatrix}
             1 & -\frac{2}{3} & 0 \\
             0 & 1 & 0 \\
             0 & 0 & 1
             \end{pmatrix}\begin{pmatrix}
             p \\
             q \\
             r 
             \end{pmatrix}
\end{equation*}
to arrive at
\begin{equation*}
\frac{d}{dt}\begin{pmatrix}
             X \\
             Y \\
             Z 
             \end{pmatrix}=\frac{1}{Z}\left\{\begin{pmatrix}
             -\frac{1}{2} & 0 & 0 \\
             0 & -2 & 0 \\
             0 & 0 & -\frac{1}{2}
             \end{pmatrix}\begin{pmatrix}
             X \\
             Y \\
             Z 
             \end{pmatrix}+\cdots\right\}.
             \end{equation*}
             In this case, the local index is $(-\frac{1}{2},-2,-\frac{1}{2})$.
             
             The $\alpha$-test,
\begin{equation}
t=t_0+\alpha T, \quad X=\alpha X_1, \quad Y=\alpha Y_1, \quad Z=\alpha Z_1, \quad \alpha \rightarrow 0,
\end{equation}
yields the following reduced system:
\begin{equation*}
\frac{d}{dt}\begin{pmatrix}
             X_1 \\
             Y_1 \\
             Z_1 
             \end{pmatrix}=\frac{1}{Z_1}\begin{pmatrix}
             -\frac{1}{2} & 0 & 0 \\
             0 & -2 & 0 \\
             0 & 0 & -\frac{1}{2}
             \end{pmatrix}\begin{pmatrix}
             X_1 \\
             Y_1 \\
             Z_1
             \end{pmatrix}.
             \end{equation*}
This system can be solved by rational functions:
$$
(X_1,Y_1,Z_1)=\left(c_2(t-2c_1),c_3(t-2c_1)^4,-\frac{1}{2}(t-2c_1) \right) \quad (c_i \in {\Bbb C}).
$$
             This suggests the possibilities that $-\frac{1}{2}$ is the residue of the formal Laurent series:
\begin{equation}
z(t)=\frac{-\frac{1}{2}}{(t-t_0)}+b_1+b_2(t-t_0)+\dots+b_n(t-t_0)^{n-1}+\cdots \quad (b_i \in {\Bbb C}),
\end{equation}
and the ratio $(\frac{-\frac{1}{2}}{-\frac{1}{2}},\frac{-2}{-\frac{1}{2}})=(1,4)$ is resonance data of the formal Laurent series of $(x(t),y(t))$, respectively. There exist meromorphic solutions with two free parameters which pass through $P$.

\section{Symmetry and holomorphy}

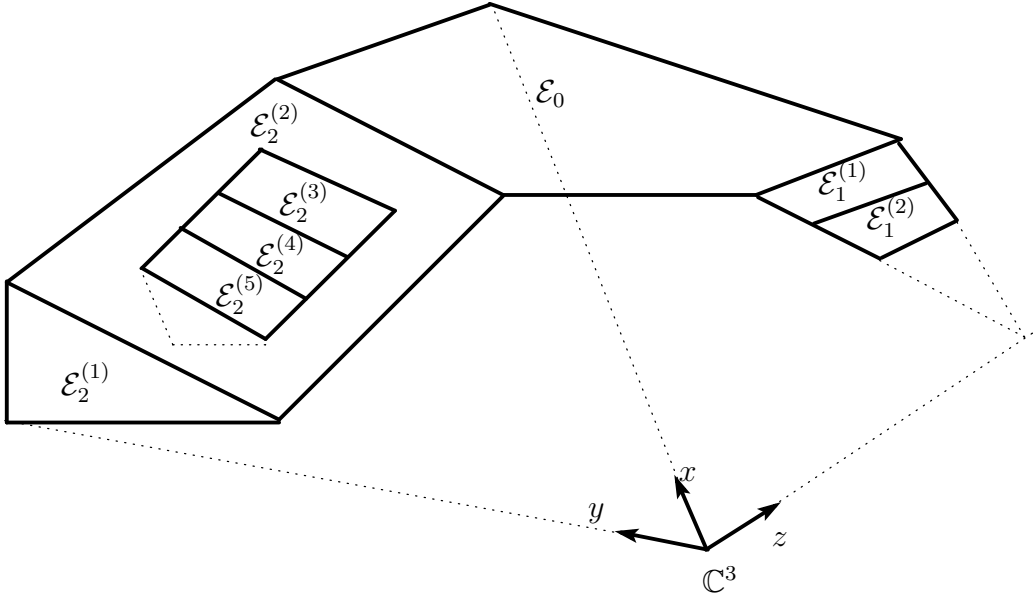
\begin{figure}[ht]
\unitlength 0.1in
\begin{picture}(53.20,29.00)(18.00,-36.50)
%
\special{pn 20}%
\special{pa 4310 755}%
\special{pa 3210 1144}%
\special{fp}%
\special{pa 3210 1144}%
\special{pa 4380 1748}%
\special{fp}%
\special{pa 4380 1748}%
\special{pa 5690 1748}%
\special{fp}%
\special{pa 5690 1748}%
\special{pa 6440 1456}%
\special{fp}%
\special{pa 6440 1456}%
\special{pa 4310 750}%
\special{fp}%
%
\special{pn 20}%
\special{pa 3200 1144}%
\special{pa 1810 2204}%
\special{fp}%
\special{pa 1810 2204}%
\special{pa 3220 2921}%
\special{fp}%
\special{pa 3220 2921}%
\special{pa 4380 1753}%
\special{fp}%
%
\special{pn 20}%
\special{pa 5990 1897}%
\special{pa 6570 1687}%
\special{fp}%
%
\special{pn 20}%
\special{pa 6340 2076}%
\special{pa 6720 1887}%
\special{fp}%
%
\special{pn 20}%
\special{pa 1810 2199}%
\special{pa 1810 2936}%
\special{fp}%
\special{pa 1810 2936}%
\special{pa 3220 2936}%
\special{fp}%
\put(45.5000,-12.9000){\makebox(0,0)[lb]{${\mathcal E}_0$}}%
\put(60.1000,-17.8200){\makebox(0,0)[lb]{${\mathcal E}_1^{(1)}$}}%
\put(62.6000,-19.6600){\makebox(0,0)[lb]{${\mathcal E}_1^{(2)}$}}%
\put(20.9000,-28.3000){\makebox(0,0)[lb]{${\mathcal E}_2^{(1)}$}}%
\put(30.7000,-14.7800){\makebox(0,0)[lb]{${\mathcal E}_2^{(2)}$}}%
\put(32.2000,-18.9000){\makebox(0,0)[lb]{${\mathcal E}_2^{(3)}$}}%
\put(28.9000,-23.7000){\makebox(0,0)[lb]{${\mathcal E}_2^{(5)}$}}%
\put(31.0000,-21.6000){\makebox(0,0)[lb]{${\mathcal E}_2^{(4)}$}}%
%
\special{pn 8}%
\special{pa 1800 2926}%
\special{pa 5440 3598}%
\special{dt 0.045}%
\special{pa 5440 3598}%
\special{pa 5439 3598}%
\special{dt 0.045}%
\special{pa 5440 3598}%
\special{pa 7120 2470}%
\special{dt 0.045}%
\special{pa 7120 2470}%
\special{pa 7119 2470}%
\special{dt 0.045}%
%
\special{pn 8}%
\special{pa 4320 750}%
\special{pa 5430 3606}%
\special{dt 0.045}%
\special{pa 5430 3606}%
\special{pa 5430 3605}%
\special{dt 0.045}%
%
\special{pn 20}%
\special{pa 3130 1510}%
\special{pa 2510 2130}%
\special{fp}%
\special{pa 2510 2130}%
\special{pa 3150 2500}%
\special{fp}%
\special{pa 3150 2500}%
\special{pa 3820 1830}%
\special{fp}%
\special{pa 3820 1830}%
\special{pa 3120 1510}%
\special{fp}%
%
\special{pn 8}%
\special{pa 2500 2130}%
\special{pa 2670 2530}%
\special{dt 0.045}%
\special{pa 2670 2530}%
\special{pa 2670 2529}%
\special{dt 0.045}%
\special{pa 2670 2530}%
\special{pa 3150 2530}%
\special{dt 0.045}%
\special{pa 3150 2530}%
\special{pa 3149 2530}%
\special{dt 0.045}%
%
\special{pn 20}%
\special{pa 2910 1740}%
\special{pa 3570 2070}%
\special{fp}%
%
\special{pn 20}%
\special{pa 2710 1920}%
\special{pa 3360 2290}%
\special{fp}%
%
\special{pn 20}%
\special{pa 5430 3600}%
\special{pa 5790 3370}%
\special{fp}%
\special{sh 1}%
\special{pa 5790 3370}%
\special{pa 5723 3389}%
\special{pa 5745 3399}%
\special{pa 5745 3423}%
\special{pa 5790 3370}%
\special{fp}%
%
\special{pn 20}%
\special{pa 5430 3590}%
\special{pa 5280 3240}%
\special{fp}%
\special{sh 1}%
\special{pa 5280 3240}%
\special{pa 5288 3309}%
\special{pa 5301 3289}%
\special{pa 5325 3293}%
\special{pa 5280 3240}%
\special{fp}%
%
\special{pn 20}%
\special{pa 5440 3600}%
\special{pa 4990 3510}%
\special{fp}%
\special{sh 1}%
\special{pa 4990 3510}%
\special{pa 5051 3543}%
\special{pa 5042 3520}%
\special{pa 5059 3503}%
\special{pa 4990 3510}%
\special{fp}%
\put(48.2000,-34.6000){\makebox(0,0)[lb]{$y$}}%
\put(52.9000,-32.4000){\makebox(0,0)[lb]{$x$}}%
\put(57.7000,-35.8000){\makebox(0,0)[lb]{$z$}}%
%
\special{pn 20}%
\special{pa 5700 1760}%
\special{pa 6330 2080}%
\special{fp}%
%
\special{pn 20}%
\special{pa 6430 1480}%
\special{pa 6730 1880}%
\special{fp}%
%
\special{pn 8}%
\special{pa 6330 2090}%
\special{pa 7080 2490}%
\special{dt 0.045}%
\special{pa 7080 2490}%
\special{pa 7079 2490}%
\special{dt 0.045}%
\special{pa 7080 2490}%
\special{pa 6730 1890}%
\special{dt 0.045}%
\special{pa 6730 1890}%
\special{pa 6730 1891}%
\special{dt 0.045}%
\put(54.1000,-38.2000){\makebox(0,0)[lb]{${\Bbb C}^3$}}%
\end{picture}%
\label{fig:HPIII1}
\caption{The part surrounding bold lines coincides with $(-{K_{\tilde{\mathcal X}}})_{red}$ (see Theorem \ref{th:5}). The symbol ${\mathcal E}_0$ denotes the proper transform of boundary divisor of ${\Bbb P}^3 \times {\Bbb C}^{*}$ by $\varphi$ and  ${\mathcal E}_i^{(j)}$ denote the exceptional divisors; ${\mathcal E}_2^{(1)}$ is isomorphic to ${\Bbb P}^2$, ${\mathcal E}_1^{(1)},{\mathcal E}_2^{(3)}$ are isomorphic to ${\Bbb F}_1$, ${\mathcal E}_1^{(2)},{\mathcal E}_2^{(4)},{\mathcal E}_2^{(5)}$ are isomorphic to ${\Bbb P}^1 \times {\Bbb P}^1$ and ${\mathcal E}_0,{\mathcal E}_2^{(2)}$ are isomorphic to the surface which can be obtained by blowing up one point in ${\Bbb P}^1 \times {\Bbb P}^1$.}
\end{figure}

\begin{Theorem}\label{th:5}
The phase space ${\mathcal X}$ for the system \eqref{eq:6} is obtained by gluing three copies of ${\mathbb C}^3 \times {\Bbb C}^{*}${\rm:\rm}
\begin{center}
${U_j} \times {\Bbb C}^{*} \cong {\mathbb C}^3 \times {\Bbb C}^{*} \ni \{(x_j,y_j,z_j,t)\},  \ \ j=0,1,2$
\end{center}
via the following birational transformations{\rm:\rm}
\begin{align}
\begin{split}
0) \ &x_0=x, \quad y_0=y, \quad z_0=z,\\
1) \ &x_1=x, \quad y_1=-\left(yz+\frac{2(2\alpha_1-1)x+2\alpha_0}{t^2}\right)z, \ z_1=\frac{1}{z},\\
2) \ &x_2=x+\frac{z}{2}, \ y_2=-((y+\frac{2}{t}x-\frac{1}{t}-\frac{4(\alpha_1-1)\alpha_1-7}{4t^2})z\\
&-\frac{2(2\alpha_1-3)t x+2(\alpha_0+1)t+4\alpha_1^2-8\alpha_1+3}{t^3}+\frac{tz+4}{4t}z^2)z,\\
&z_2=\frac{1}{z}.
\end{split}
\end{align}
These transition functions satisfy the condition{\rm:\rm}
\begin{equation*}
dx_i \wedge dy_i \wedge dz_i=dx \wedge dy \wedge dz \quad (i=1,2).
\end{equation*}
\end{Theorem}
We note that ${\Bbb C}^{*}={\Bbb C}-\{0\}$.

\begin{Theorem}\label{th:6}
After nine successive blowing-ups in ${\Bbb P}^3 \times {\Bbb C}^{*}$, we obtain the smooth projective $4$-fold $\tilde{\mathcal X}$ and a morphism $\varphi:\tilde{\mathcal X} \rightarrow {\Bbb P}^3 \times {\Bbb C}^{*}$. Its canonical divisor $K_{\tilde{\mathcal X}}$ is given by
\begin{align}
\begin{split}
K_{\tilde{\mathcal X}}&=-4{\mathcal E}_0 -2{\mathcal E}_1^{(1)}-{\mathcal E}_1^{(2)}-2{\mathcal E}_2^{(1)}-5{\mathcal E}_2^{(2)}-3{\mathcal E}_2^{(3)}-2{\mathcal E}_2^{(4)}-{\mathcal E}_2^{(5)},
\end{split}
\end{align}
where the symbol ${\mathcal E}_0$ denotes the proper transform of boundary divisor of ${\Bbb P}^3 \times {\Bbb C}^{*}$ by $\varphi$ and  ${\mathcal E}_i^{(j)}$ denote the exceptional divisors; ${\mathcal E}_2^{(1)}$ is isomorphic to ${\Bbb P}^2$, ${\mathcal E}_1^{(1)},{\mathcal E}_2^{(3)}$ are isomorphic to ${\Bbb F}_1$, ${\mathcal E}_1^{(2)},{\mathcal E}_2^{(4)},{\mathcal E}_2^{(5)}$ are isomorphic to ${\Bbb P}^1 \times {\Bbb P}^1$ and ${\mathcal E}_0,{\mathcal E}_2^{(2)}$ are isomorphic to the surface which can be obtained by blowing up one point in ${\Bbb P}^1 \times {\Bbb P}^1$ \rm{ (see Figure 1) \rm}. Moreover, $\tilde{\mathcal X}-(-{K_{\tilde{\mathcal X}}})_{red}$ satisfies
\begin{equation}
\tilde{\mathcal X}-(-{K_{\tilde{\mathcal X}}})_{red}={\mathcal X}.
\end{equation}
\end{Theorem}

{\bf Proof of Theorems \ref{th:5} and \ref{th:6}.}

By the following steps, we can resolve the accessible singular point
$$
Q:=\{(X_2,Y_2,Z_2)=(0,0,0)\}.
$$

{\bf Step 1}: We blow up at $Q${\rm : \rm}
$$
p^{(1)}=\frac{X_2}{Z_2}, \quad q^{(1)}=Z_2, \quad r^{(1)}=\frac{Y_2}{Z_2}.
$$

{\bf Step 2}: We blow up along the curve $\{(p^{(1)},q^{(1)},r^{(1)})=(p^{(1)},0,0)\}${\rm : \rm}
$$
p^{(2)}=p^{(1)}, \quad q^{(2)}=\frac{q^{(1)}}{r^{(1)}}, \quad r^{(2)}=r^{(1)}.
$$

{\bf Step 3}: We make the change of variables{\rm : \rm}
$$
p^{(3)}=p^{(2)}, \quad q^{(3)}=\frac{1}{q^{(2)}}, \quad r^{(3)}=r^{(2)}.
$$

{\bf Step 4}: We blow up at the point $\{(p^{(3)},q^{(3)},r^{(3)})=(-\frac{1}{2},-\frac{1}{4},0) \}${\rm : \rm}
$$
p^{(4)}=\frac{p^{(3)}+\frac{1}{2}}{r^{(3)}}, \quad q^{(4)}=\frac{q^{(3)}+\frac{1}{4}}{r^{(3)}}, \quad r^{(4)}=r^{(3)}.
$$

{\bf Step 5}: We blow up along the curve $\{(p^{(4)},q^{(4)},r^{(4)})|q^{(4)}=r^{(4)}=0 \}${\rm : \rm}
$$
p^{(5)}=p^{(4)}, \quad q^{(5)}=\frac{q^{(4)}}{r^{(4)}}, \quad r^{(5)}=r^{(4)}.
$$

{\bf Step 6}: We blow up along the curve $\{(p^{(5)},q^{(5)},r^{(5)})|q^{(5)}-\frac{5-12\alpha_1+4\alpha_1^2+4t-8t p^{(5)}}{4t^2}\\
=r^{(5)}=0 \}${\rm : \rm}
$$
p^{(6)}=p^{(5)}, \quad q^{(6)}=\frac{q^{(5)}-\frac{5-12\alpha_1+4\alpha_1^2+4t-8t p^{(5)}}{4t^2}}{r^{(5)}}, \quad r^{(6)}=r^{(5)}.
$$

{\bf Step 7}: We blow up along the curve $\{(p^{(6)},q^{(6)},r^{(6)})|\\
q^{(6)}-\frac{3-8\alpha_1+4\alpha_1^2+2t+2\alpha_0 t-6t p^{(6)}+4\alpha_1 t p^{(6)} }{t^3}=r^{(6)}=0 \}${\rm : \rm}
$$
p^{(7)}=p^{(6)}, \quad q^{(7)}=\frac{q^{(6)}-\frac{3-8\alpha_1+4\alpha_1^2+2t+2\alpha_0 t-6t p^{(6)}+4\alpha_1 t p^{(6)} }{t^3}}{r^{(6)}}, \quad r^{(7)}=r^{(6)}.
$$
We have resolved the accessible singular point $Q$. By choosing a new coordinate system as
$$
(x_2,y_2,z_2)=(p^{(7)},-q^{(7)},r^{(7)}),
$$
we can obtain the coordinate system $(x_2,y_2,z_2)$ in the description of $\mathcal X$ given in Theorem \ref{th:5}.

The proof on the accessible singular point $P_2$ is same.

Now we complete the proof of Theorems \ref{th:5} and \ref{th:6}. \qed

\begin{Proposition}\label{pro:3}
Let us consider a system of first order ordinary differential equations in the polynomial class\rm{:\rm}
\begin{equation*}
\frac{dx}{dt}=f_1(x,y,z), \quad \frac{dy}{dt}=f_2(x,y,z), \quad \frac{dz}{dt}=f_3(x,y,z).
\end{equation*}
We assume that

$(A1)$ $deg(f_i)=3$ with respect to $x,y,z$.

$(A2)$ The right-hand side of this system becomes again a polynomial in each coordinate system $(x_i,y_i,z_i) \ (i=1,2)$.

\noindent
Then such a system coincides with the system \eqref{eq:6}.
\end{Proposition}

\begin{Theorem}\label{th:2}
The system \eqref{eq:6} is invariant under the following transformations\rm{: \rm}
\begin{align*}
\begin{split}
&s_0:(x,y,z,t;\alpha_0,\alpha_1) \rightarrow \left(x,y,z+\frac{2(2\alpha_1-1)x+2\alpha_0}{t^2y},t;-\alpha_0,1-\alpha_1 \right),\\
&s_1:\left(x,y,z,t;\alpha_0,\frac{3}{2} \right) \rightarrow \left(x+2G,y+2G-4G^2,z-4G,t;-2-\alpha_0,\frac{3}{2} \right),\\
&\pi: \left(x,y,z,t;\alpha_0,\frac{3}{2} \right)\\
&\rightarrow \left(x+\frac{z}{2},-y-\frac{2x}{t}-\frac{z^2}{4}-\frac{z}{t}+\frac{(2\alpha_1-1)(2\alpha_1-3)}{4t^2}+\frac{1}{t},-z-\frac{2}{t},t;-\alpha_0-1,\frac{1}{2} \right),
\end{split}
\end{align*}
where the symbol $G$ denotes
\begin{equation}
G:=\frac{t\left\{2(2\alpha_1-3)x+(2\alpha_1-3)z+2(\alpha_0+1) \right\}+(2\alpha_1-1)(2\alpha_1-3)}{t\left\{t^2(4y+z^2)+4t(2x+z-1)-(2\alpha_1-1)(2\alpha_1-5) \right\}}.
\end{equation}
\end{Theorem}
We note that the transformations $s_1$ and $\pi$ become its B{\"a}cklund transformations for the system \eqref{eq:6} only if $\alpha_1=\frac{3}{2}$, respectively.

\begin{Proposition}\label{pro:2}
The transformation $s_0s_1$ acts on $\alpha_0$ as
\begin{equation}
s_0s_1(\alpha_0)=\alpha_0-2.
\end{equation}
\end{Proposition}

\section{Special solutions}
In this section, we study a solution of the system \eqref{eq:6} which is written by the use of classical transcendental functions.

By setting $(\alpha_0,\alpha_1)=\left(0,\frac{3}{2} \right)$, we obtain a solution of the form:
\begin{align}
x=0, \quad y=0, \quad \frac{dz}{dt}=-\frac{1}{2}z^2-\frac{2}{t}z+\frac{2}{t}.
\end{align}
By making the change of variables
\begin{equation*}
z=2\frac{d}{dt}\rm{log \rm}Z,
\end{equation*}
we obtain the second-order linear equation:
\begin{equation}
t\frac{d^2Z}{dt^2}+2\frac{dZ}{dt}-Z=0.
\end{equation}
This equation can be solved by
\begin{equation}
Z(t)=\sum_{k=0}^{\infty} \frac{1}{(2)_k}\frac{t^k}{k!}.
\end{equation}

\end{document}